\documentclass{amsproc}%
\usepackage{graphicx}
\usepackage{amscd}
\usepackage{amsmath}
\usepackage{amsfonts}
\usepackage{amssymb}%
\theoremstyle{plain}

\numberwithin{equation}{section}

\begin{document}
\title[On Lie nilpotent rings and Cohen's Theorem ]{On Lie nilpotent rings and Cohen's Theorem}
\author{Jen\H{o} Szigeti}
\address{Institute of Mathematics, University of Miskolc, Miskolc, Hungary 3515}
\email{jeno.szigeti@uni-miskolc.hu}
\author{Leon van Wyk}
\address{Department of Mathematical Sciences, Stellenbosch University\\
P/Bag X1, Matieland 7602, Stellenbosch, South Africa }
\email{LvW@sun.ac.za}
\thanks{The first author was supported by OTKA K-101515 of Hungary and by the
TAMOP-4.2.1.B-10/2/KONV-2010-0001 project with support by the European Union,
co-financed by the European Social Fund.}
\thanks{The second author was supported by the National Research Foundation of South
Africa under Grant No.~UID 72375. Any opinion, findings and conclusions or
recommendations expressed in this material are those of the authors and
therefore the National Research Foundation does not accept any liability in
regard thereto.}
\thanks{The authors thank P. N. Anh, L. Marki and J. H. Meyer for fruitful consultations.}

\begin{abstract}
We study certain (two-sided) nil ideals and nilpotent ideals in a Lie
nilpotent ring $R$. Our results lead us to showing that the prime radical
$\mathrm{rad}(R)$ of $R$ comprises the nilpotent elements of $R$, and that if
$L$ is a left ideal of $R$, then $L+\mathrm{rad}(R)$ is a two-sided ideal of
$R$. This in turn leads to a Lie nilpotent version of Cohen's theorem, namely
if $R$ is a Lie nilpotent ring and every prime (two-sided) ideal of $R$ is
finitely generated as a left ideal, then every left ideal of $R$ containing
the prime radical of $R$ is finitely generated (as a left ideal). For an
arbitrary ring~$R$ with identity we also consider its so-called $n$-th Lie
center $Z_{n}(R),\ n\geq1$, which is a Lie nilpotent ring of index $n$. We
prove that if $C$ is a commutative submonoid of the multiplicative monoid
of~$R$, then the subring $\langle\mathrm{Z}_{n}(R)\cup C\rangle$ of~$R$
generated by the subset $\mathrm{Z}_{n}(R)\cup C$ of~$R$ is also Lie nilpotent
of index $n.$

\end{abstract}
\maketitle

\noindent1. INTRODUCTION

\bigskip

\noindent One of the inspirations for this paper was a search for a Lie
nilpotent version of a well known theorem in commutative algebra, namely
Cohen's Theorem (see [Co]), which states that a commutative ring $R$ is
Noetherian (i.e.~every ideal of $R$ is finitely generated) if every prime
ideal of $R$ is finitely generated.

Cohen's Theorem and Kaplansky's Theorem (see [Ka]) are quite often mentioned
in the same breath. The latter states that a commutative Noetherian ring $R$
is a principal ideal ring (i.e.~every left ideal of $R$ is principal and every
right ideal of~$R$ is principal) if and only if every maximal ideal of $R$ is
principal. In fact, in [Re], an excellent paper which contains a thorough
survey (including an extensive list of references) of quite a number of
versions of Cohen's Theorem in various contexts, mention is made of the
combined Kaplansky-Cohen Theorem, which states that a commutative ring $R$ is
a principal ideal ring if and only if every prime ideal of $R$ is principal.

Moreover, Cohen's Theorem and the Kaplansky-Cohen Theorem are strengthened in
[Re], and some generalizations of Cohen's Theorem, for example, Koh's (see
[Ko]) and Chandran's (see [Ch]), are implied by results in [Re]. Other papers
on Cohen's Theorem, some of of which contain module versions of Cohen's
Theorem, include [H], [Jo], [L], [Na] and [Ni].

In Section 2 we deal with certain products in a Lie nilpotent ring $R$ of
index $n\geq2$. We shall make use of classical results due to Jennings (see
[Je]), and in order to ease readability, we provide short self-contained
proofs of these theorems. A basic example of a Lie nilpotent ring $R$ of index
$n\geq2$ is also presented.

In Section 3 we show, amongst others, that the prime radical $\mathrm{rad}(R)$
of a Lie nilpotent ring $R$ of index $n\geq2$ comprises the nilpotent elements
of $R$, and that if $L$ is a left ideal of $R$, then $L+\mathrm{rad}(R)$ is a
two-sided ideal of $R$. This in turn leads to Theorem 3.3: if $R$ is a Lie
nilpotent ring and every prime (two-sided) ideal of $R$ is finitely generated
as a left ideal, then every left ideal of $R$ containing the prime radical of
$R$ is finitely generated (as a left ideal). In this regard we mention (see
[Re] and [Kr]) that G. Michler and L. Small proved independently that a left
fully bounded ring (such as a polynomial identity ring) in which every prime
ideal is finitely generated as a left ideal is left noetherian. Therefore,
although Theorem 3.3 is not new, the particular techniques employed in the
present paper in the study of Lie nilpotent rings are important in their own right.

The authors note that they use the notation $\mathrm{rad}(R)$ for the prime
(or lower nil) radical and $\mathrm{J}(R)$ for the Jacobson radical of $R$.

Another inspiration for this paper was Lemma 2.1 of [Sz], which plays a
crucial role in the development of the Lie nilpotent determinant theory in
[Sz] and [SzvW]. For an arbitrary ring $R$ with identity we consider in
Section 4 its so-called $n$-th Lie center $\mathrm{Z}_{n}(R),\ n\geq1$, which
is a Lie nilpotent ring of index $n$. We prove the following broad
generalization of the mentioned lemma: if $C$ is a commutative submonoid of
the multiplicative monoid of $R$, then the subring $\langle\mathrm{Z}%
_{n}(R)\cup C\rangle$ of $R$ generated by the subset $\mathrm{Z}_{n}(R)\cup C$
of $R$ is also Lie nilpotent of index $n$.

\bigskip

\noindent2. PRODUCTS IN LIE NILPOTENT RINGS

\bigskip

Let $R$ be a ring, and let $[x,y]=xy-yx$ denote the additive commutator of the
elements $x,y\in R$. It is well known that $(R,+,[$ $,$ $])$ is a Lie ring and
the following identities hold:%
\[
\lbrack y,x]=-[x,y],
\]%
\[
\lbrack\lbrack x,y],z]+[[y,z],x]+[[z,x],y]=0\text{ (Jacobian identity)},
\]%
\[
\lbrack uv,x]=[u,vx]+[v,xu]=u[v,x]+[u,x]v.
\]
The use of $[x,uv]=u[x,v]+[x,u]v$ and $[uv,x]=u[v,x]+[u,x]v$\ gives that%
\[
\lbrack\lbrack a,by],x]=[b[a,y]+[a,b]y,x]=[b[a,y],x]+[[a,b]y,x]=
\]%
\[
b[[a,y],x]+[b,x][a,y]+[a,b][y,x]+[[a,b],x]y.
\]
For a sequence $x_{1},x_{2},\ldots,x_{n}$ of elements in $R$ we use the
notation $[x_{1},x_{2},\ldots,x_{n}]_{n}^{\ast}$\ for the \textit{left normed
commutator (Lie-)product}:%
\[
\lbrack x_{1}]_{1}^{\ast}=x_{1}\text{ and }[x_{1},x_{2},\ldots,x_{n}%
]_{n}^{\ast}=[\ldots\lbrack\lbrack x_{1},x_{2}],x_{3}],\ldots,x_{n}].
\]
Clearly, we have%
\[
\lbrack x_{1},x_{2},\ldots,x_{n},x_{n+1}]_{n+1}^{\ast}=[[x_{1},x_{2}%
,\ldots,x_{n}]_{n}^{\ast},x_{n+1}]=[[x_{1},x_{2}],x_{3},\ldots,x_{n}%
,x_{n+1}]_{n}^{\ast}.
\]
A ring $R$ is called \textit{Lie nilpotent of index }$n$ (or \textit{having
property} $\mathrm{L}_{n}$) if%
\[
\lbrack x_{1},x_{2},\ldots,x_{n},x_{n+1}]_{n+1}^{\ast}=0
\]
is a polynomial identity on $R$. If $R$ has property $\mathrm{L}_{n}$, then
$[x_{1},x_{2},\ldots,x_{n}]_{n}^{\ast}\in\mathrm{Z}(R)$ is central for all
$x_{1},x_{2},\ldots,x_{n}\in R$.

Let $k_{0}=0$ and $1\leq k_{1},\ldots,k_{n},k_{n+1}\leq m$ be integers such
that $k_{1}+\cdots+k_{n}+k_{n+1}=m$ and let $K$ be a field. The pair $(i,j)$
of integers \textit{satisfies }$(\ast)$ if%
\[
(\ast)\text{ \ \ \ \ \ \ }k_{0}+k_{1}+\cdots+k_{t-1}<i\leq k_{0}+k_{1}%
+\cdots+k_{t}<j\leq m
\]
for some (unique) index $1\leq t\leq n$. One of the basic examples of a
$K$-algebra satisfying $\mathrm{L}_{n}$ is the $K$-subalgebra%
\[
R=R_{m}(k_{1},\ldots,k_{n},k_{n+1})=\left\{
{\textstyle\sum}
a_{i,j}E_{i,j}\mid a_{i,j}\in K\text{ and }(i,j)\text{ satisfies }%
(\ast)\right\}
\]
of block upper triangular $m\times m$\ matrices of the full matrix algebra
$\mathrm{M}_{m}(K)$, where $E_{i,j}$ is the standard matrix unit with $1$ in
the $(i,j)$ position. Clearly, the ordinary nilpotency $R^{n+1}=\{0\}\neq
R^{n}$ implies $\mathrm{L}_{n}$, and the addition of the center (scalar
matrices) yields a unitary subring $R+KI_{m}$ of $\mathrm{M}_{m}(K)$ also with
$\mathrm{L}_{n}$. The $K$-dimension of $R+KI_{m}$ is%
\[
\dim_{K}(R)=1+\frac{1}{2}(m^{2}-k_{1}^{2}-\cdots-k_{n}^{2}-k_{n+1}^{2}).
\]
\noindent\textbf{Conjecture.} \textit{If }$S$\textit{ is a (unitary) }%
$K$\textit{-subalgebra of }$\mathrm{M}_{m}(K)$\textit{ with }$L_{n}$\textit{
and }$n+1\leq m$\textit{, then}%
\[
\dim_{K}(S)\leq1+\frac{1}{2}(m^{2}-k_{1}^{2}-\cdots-k_{n}^{2}-k_{n+1}^{2})
\]
\textit{for some integers }$1\leq k_{1},\ldots,k_{n},k_{n+1}\leq m$\textit{
with }$k_{1}+\cdots+k_{n}+k_{n+1}=m$\textit{.}

\bigskip

For $n=1$ our conjecture becomes a classical theorem of Schur about the
maximal dimension of a commutative subalgebra of $\mathrm{M}_{m}(K)$.

We call a ring $R$ \textit{Lie nilpotent} if it is Lie nilpotent of index $n$
for some $n\geq1$ (see, for example, [D], [DF],[Je] and [Ro]). Notice that
$\mathrm{M}_{2}(K)$ is not Lie nilpotent, i.e. it does not satisfy
$\mathrm{L}_{n}$ for any $n\geq1$.

\bigskip

\noindent\textbf{2.1.Proposition.} \textit{Let }$R$\textit{ be a ring with
}$\mathrm{L}_{2}$\textit{ and }$a,b,c,a_{0},a_{1},\ldots,a_{k}\in R$\textit{.}

\begin{enumerate}
\item $[a,b][a,c]=0$\textit{.}

\item \textit{If }$ab=0$\textit{, then }$bxbya=0$\textit{ for all }$x,y\in
R$\textit{ (i.e. }$bRbRa=\{0\}$\textit{).}

\item \textit{If }$a_{0}a_{1}\cdots a_{k}=0$\textit{,} \textit{then}%
\[
a_{1}x_{1}a_{1}y_{1}a_{2}x_{2}a_{2}y_{2}\cdots a_{k}x_{k}a_{k}y_{k}a_{0}=0
\]
\textit{for all }$x_{i},y_{i}\in R$\textit{, }$1\leq i\leq k$\textit{.}
\end{enumerate}

\bigskip

\noindent\textbf{Proof.}

(1): Take $x=c$ and $y=a$ in%
\[
\lbrack\lbrack a,by],x]=b[[a,y],x]+[b,x][a,y]+[a,b][y,x]+[[a,b],x]y.
\]

(2): $ab=0$ and the $\mathrm{L}_{2}$ property of $R$ imply that $bya=[by,a]$
is central, whence%
\[
bxbya=(bya)bx=0
\]
follows.

(3): In order to see the validity of the implication for $k=1$, take $a=a_{0}$
and $b=a_{1}$ in part (2). In the next step of the induction we assume that
our statement holds for some $k\geq1$ and consider the product $a_{0}%
a_{1}\cdots a_{k}a_{k+1}=0$ in $R$. Using the induction hypothesis for
$(a_{0}a_{1})a_{2}\cdots a_{k}a_{k+1}=0$, we obtain that%
\[
a_{2}x_{2}a_{2}y_{2}a_{3}\cdots a_{k}x_{k}a_{k}y_{k}a_{k+1}x_{k+1}%
a_{k+1}y_{k+1}(a_{0}a_{1})=0
\]
for all $x_{i},y_{i}\in R$, $2\leq i\leq k+1$. Now the choice of%
\[
a=a_{2}x_{2}a_{2}y_{2}a_{3}\cdots a_{k}x_{k}a_{k}y_{k}a_{k+1}x_{k+1}%
a_{k+1}y_{k+1}a_{0}\text{ and }b=a_{1}%
\]
in part (2) gives that%
\[
0=bx_{1}by_{1}a=a_{1}x_{1}a_{1}y_{1}a_{2}x_{2}a_{2}y_{2}a_{3}\cdots a_{k}%
x_{k}a_{k}y_{k}a_{k+1}x_{k+1}a_{k+1}y_{k+1}a_{0}%
\]
for all $x_{i},y_{i}\in R$, $1\leq i\leq k+1$. It follows that our statement
is valid for $k+1$. $\square$

\bigskip

\noindent\textbf{2.2.Theorem ([Je]).} \textit{Let }$n\geq3$\textit{ be an
integer and }$R$\textit{ be a ring with }$\mathrm{L}_{n}$\textit{. Then}%
\[
\lbrack x_{1},x_{2},\ldots,x_{n}]_{n}^{\ast}\cdot\lbrack y_{1},y_{2}%
,\ldots,y_{n}]_{n}^{\ast}=0
\]
\textit{for all }$x_{i},y_{i}\in R$\textit{, }$1\leq i\leq n$\textit{. Thus
the two sided ideal}%
\[
N\!=\!R\{[x_{1},x_{2},\ldots,x_{n}]_{n}^{\ast}\!\mid\!x_{i}\in R,1\leq i\leq
n\}\!=\!\{[x_{1},x_{2},\ldots,x_{n}]_{n}^{\ast}\!\mid\!x_{i}\in R,1\leq i\leq
n\}R
\]
\textit{generated by the (central) elements }$[x_{1},x_{2},\ldots,x_{n}%
]_{n}^{\ast}$\textit{ is nilpotent with }$N^{2}=\{0\}$\textit{.}

\bigskip

\noindent\textbf{Proof.} Take%
\[
x=[x_{1},x_{2},\ldots,x_{n-2}]_{n-2}^{\ast},y=x_{n-1}\text{ and }%
z=[y_{1},y_{2},\ldots,y_{n-1}]_{n-1}^{\ast}%
\]
in the Jacobi identity%
\[
\lbrack\lbrack x,y],z]+[[y,z],x]+[[z,x],y]=0.
\]
Since%
\[
\lbrack\lbrack y,z],x]\!=\!-\![[z,y],x]\!=\!-\![[[y_{1},y_{2},\ldots
,y_{n-1}]_{n-1}^{\ast},y],x]\!=\!-\![y_{1},y_{2},\ldots,y_{n-1},y,x]_{n+1}%
^{\ast}\!=\!0
\]
and%
\[
\lbrack\lbrack z,x],y]=[[[y_{1},y_{2},\ldots,y_{n-1}]_{n-1}^{\ast
},x],y]=[y_{1},y_{2},\ldots,y_{n-1},x,y]_{n+1}^{\ast}=0
\]
are consequences of the $\mathrm{L}_{n}$ property, we obtain that%
\[
0=[[x,y],z]=[[[x_{1},x_{2},\ldots,x_{n-2}]_{n-2}^{\ast},x_{n-1}],[y_{1}%
,y_{2},\ldots,y_{n-1}]_{n-1}^{\ast}]=
\]%
\[
\lbrack\lbrack x_{1},x_{2},\ldots,x_{n-1}]_{n-1}^{\ast},[y_{1},y_{2}%
,\ldots,y_{n-1}]_{n-1}^{\ast}].
\]
Now take%
\[
a=[y_{1},y_{2},\ldots,y_{n-1}]_{n-1}^{\ast},b=[x_{1},x_{2},\ldots
,x_{n-1}]_{n-1}^{\ast},x=x_{n}\text{ and }y=y_{n}%
\]
in%
\[
\lbrack\lbrack a,by],x]=b[[a,y],x]+[b,x][a,y]+[a,b][y,x]+[[a,b],x]y.
\]
Since%
\[
\lbrack\lbrack a,by],x]=[[[y_{1},y_{2},\ldots,y_{n-1}]_{n-1}^{\ast
},by],x]=[y_{1},y_{2},\ldots,y_{n-1},by,x]_{n+1}^{\ast}=0,
\]%
\[
\lbrack\lbrack a,y],x]=[[[y_{1},y_{2},\ldots,y_{n-1}]_{n-1}^{\ast
},y],x]=[y_{1},y_{2},\ldots,y_{n-1},y,x]_{n+1}^{\ast}=0,
\]%
\[
\lbrack\lbrack a,b],x]=[[[y_{1},y_{2},\ldots,y_{n-1}]_{n-1}^{\ast
},b],x]=[y_{1},y_{2},\ldots,y_{n-1},b,x]_{n+1}^{\ast}=0
\]
are consequences of the $\mathrm{L}_{n}$ property and%
\[
\lbrack a,b]=[[y_{1},y_{2},\ldots,y_{n-1}]_{n-1}^{\ast},[x_{1},x_{2}%
,\ldots,x_{n-1}]_{n-1}^{\ast}]=0,
\]
we obtain that%
\[
0=[b,x][a,y]=[[x_{1},x_{2},\ldots,x_{n-1}]_{n-1}^{\ast},x_{n}][[y_{1}%
,y_{2},\ldots,y_{n-1}]_{n-1}^{\ast},y_{n}]=
\]%
\[
\lbrack x_{1},x_{2},\ldots,x_{n}]_{n}^{\ast}[y_{1},y_{2},\ldots,y_{n}%
]_{n}^{\ast}.\text{ }\square
\]

\noindent\textbf{2.3.Remark.} The $m$-generated ($m\geq4$) Grassmann algebra
\[
E^{(m)}=K\left\langle v_{1},\ldots,v_{m}\mid v_{i}v_{j}+v_{j}v_{i}=0\text{ for
all }1\leq i\leq j\leq m\right\rangle
\]
over a field $K$ (with $4\neq0$) has property $\mathrm{L}_{2}$ and%
\[
\lbrack v_{1},v_{2}]_{2}^{\ast}\cdot\lbrack v_{3},v_{4}]_{2}^{\ast}%
=[v_{1},v_{2}]\cdot\lbrack v_{3},v_{4}]=4v_{1}v_{2}v_{3}v_{4}\neq0
\]
shows that Theorem 2.2 is not valid for $n=2$. The fact that $a$ occurs twice
in $[a,b][a,c]$ is essential in part (1) of Proposition 2.1.

\bigskip

\noindent\textbf{2.4.Corollary ([Je]).} \textit{Let }$n\geq2$\textit{ be an
integer and }$R$\textit{ be a ring with }$\mathrm{L}_{n}$\textit{.}

\begin{enumerate}
\item \textit{The ideal}%
\[
I(2)=R[R,R]R=\{%
{\textstyle\sum\nolimits_{1\leq i\leq t}}
r_{i}[a_{i},b_{i}]s_{i}\mid r_{i},a_{i},b_{i},s_{i}\in R,1\leq i\leq
t\}\vartriangleleft R
\]
\textit{is nil.}

\item \textit{The ideal}%
\[
I(3)=R[[R,R],R]R=\{%
{\textstyle\sum\nolimits_{1\leq i\leq t}}
r_{i}[[a_{i},b_{i}],c_{i}]s_{i}\mid r_{i},a_{i},b_{i},c_{i},s_{i}\in R,1\leq
i\leq t\}\vartriangleleft R
\]
\textit{is nilpotent of index }$2^{n-2}$\textit{.}
\end{enumerate}

\bigskip

\noindent\textbf{Proof.} In both cases we use the ideal $N\vartriangleleft R$
generated by the (central) elements $[x_{1},x_{2},\ldots,x_{n}]_{n}^{\ast}$
and apply an induction with respect to $n$.

(1): If $n=2$, then part (1) of Proposition 2.1 gives that $[a_{i},b_{i}%
]^{2}=0$ for each $1\leq i\leq t$. Thus%
\[
\left(
{\textstyle\sum\nolimits_{1\leq i\leq t}}
r_{i}[a_{i},b_{i}]s_{i}\right)  ^{t+1}=0
\]
follows from the centrality of the $[a_{i},b_{i}]$'s. Now assume that (1) is
valid in any ring with $\mathrm{L}_{n-1}$ (for some $n\geq3$) and consider an
element%
\[%
{\textstyle\sum\nolimits_{1\leq i\leq t}}
r_{i}[a_{i},b_{i}]s_{i}\in R[R,R]R
\]
in a ring $R$ with $\mathrm{L}_{n}$. Since the factor ring $\overline{R}=R/N$
satisfies $\mathrm{L}_{n-1}$, the induction hypothesis gives that for some
$k\geq1$ we have%
\[
\left(
{\textstyle\sum\nolimits_{1\leq i\leq t}}
\overline{r_{i}}[\overline{a_{i}},\overline{b_{i}}]\overline{s_{i}}\right)
^{k}=\overline{0}%
\]
in $\overline{R}$. Since $N^{2}=\{0\}$ by Theorem 2.2, we have%
\[
\left(
{\textstyle\sum\nolimits_{1\leq i\leq t}}
r_{i}[a_{i},b_{i}]s_{i}\right)  ^{2k}=0
\]
in $R$.

(2): If $n=2$, then $R[[R,R],R]R=\{0\}$. Now assume that (2) is valid in any
ring with $\mathrm{L}_{n-1}$ ($n\geq3$) and consider a sequence of elements
$[[a_{i},b_{i}],c_{i}]$, $1\leq i\leq q=2^{n-2}$ in $R[[R,R],R]R$, where $R$
is a ring with $\mathrm{L}_{n}$. Since the factor ring $\overline{R}=R/N$
satisfies $\mathrm{L}_{n-1}$, the induction hypothesis gives that for
$p=2^{n-3}$%
\[
\lbrack\lbrack\overline{a_{1}},\overline{b_{1}}],\overline{c_{1}}]\overline
{R}[[\overline{a_{2}},\overline{b_{2}}],\overline{c_{2}}]\overline{R}%
\cdots\overline{R}[[\overline{a_{p}},\overline{b_{p}}],\overline{c_{p}%
}]=\{\overline{0}\}
\]
holds in $\overline{R}$. Thus%
\[
\lbrack\lbrack a_{1},b_{1}],c_{1}]R[[a_{2},b_{2}],c_{2}]R\cdots R[[a_{p}%
,b_{p}],c_{p}]\subseteq N
\]
and similarly%
\[
\lbrack\lbrack a_{p+1},b_{p+1}],c_{p+1}]R[[a_{p+2},b_{p+2}],c_{p+2}]R\cdots
R[[a_{p+p},b_{p+p}],c_{p+p}]\subseteq N
\]
for all $a_{i},b_{i},c_{i}\in R$, $1\leq i\leq2p=q$. It follows that%
\[
\lbrack\lbrack a_{1},b_{1}],c_{1}]R\cdots R[[a_{p},b_{p}],c_{p}]R[[a_{p+1}%
,b_{p+1}],c_{p+1}]R\cdots R[[a_{p+p},b_{p+p}],c_{p+p}]\subseteq N^{2}.
\]
In view of Theorem 2.2, we obtain that%
\[
(R[[R,R],R]R)^{q}\subseteq N^{2}=\{0\}.\text{ }\square
\]

\bigskip

\noindent\textbf{2.5.Theorem.} \textit{Let }$n\geq2$\textit{ be an integer and
}$R$\textit{ be a ring with }$\mathrm{L}_{n}$\textit{. If }$a,b\in R$\textit{,
}$ab=0$\textit{ and }$q=q(n)=2^{n-2}$, \textit{then}%
\[
bx_{1}by_{1}az_{1}bx_{2}by_{2}az_{2}\cdots z_{q-1}bx_{q}by_{q}%
a=0\text{\textit{, i.e. }}\underset{1.}{\underbrace{bRbRa}}R\underset
{2.}{\underbrace{bRbRa}}R\ldots R\underset{q.}{\underbrace{bRbRa}}=\{0\}
\]
\textit{for all }$x_{i},y_{i},z_{i}\in R$\textit{, }$1\leq i\leq q$\textit{
(}$z_{q}=1$\textit{).}

\bigskip

\noindent\textbf{Proof.} If $n=2$, then $q(2)=1$ and part (2) of Proposition
2.1 gives the result. Now assume (by induction) that our statement is valid in
any ring with $\mathrm{L}_{n-1}$ (for some $n\geq3$) and consider the elements
$a,b\in R$ with $ab=0$ in a ring $R$ with $\mathrm{L}_{n}$. In view of Theorem
2.2, we have $N^{2}=\{0\}$ for the ideal $N\vartriangleleft R$ generated by
the (central) elements $[x_{1},x_{2},\ldots,x_{n}]_{n}^{\ast}$. Since the
factor ring $R/N$ satisfies $\mathrm{L}_{n-1}$, the induction hypothesis gives
that%
\[
(b+N)(x_{1}+N)(b+N)(y_{1}+N)(a+N)(z_{1}+N)\cdots
\]%
\[
\cdots(z_{q(n-1)-1}+N)(b+N)(x_{q(n-1)}+N)(b+N)(y_{q(n-1)}+N)(a+N)=0+N
\]
in $R/N$. Thus%
\[
bx_{1}by_{1}az_{1}\cdots z_{q(n-1)-1}bx_{q(n-1)}by_{q(n-1)}a\in N
\]
and similarly%
\[
bx_{q(\!n\!-\!1\!)\!+\!1}by_{q(\!n\!-\!1\!)\!+\!1}az_{q(\!n\!-\!1\!)\!+\!1}%
\!\cdots\!z_{q(\!n\!-\!1\!)\!+\!q(\!n\!-\!1\!)\!-\!1}%
bx_{q(\!n\!-\!1\!)\!+\!q(\!n\!-\!1\!)}by_{q(\!n\!-\!1\!)\!+\!q(\!n\!-\!1\!)}%
a\!\in\!N
\]
for all $x_{i},y_{i}\in R$, $1\leq i\leq2q(n-1)$ and $z_{i}\in R$, $1\leq
i\leq2q(n-1)-1$, $i\neq q(n-1)$. Now for any $z_{q(n-1)}\in R$ we have%
\[
(bx_{1}by_{1}az_{1}\!\cdots\!z_{q(n-1)-1}bx_{q(n-1)}by_{q(n-1)}a)z_{q(n-1)}%
(bx_{q(n-1)+1}by_{q(n-1)+1}az_{q(n-1)+1}\!\cdots
\]%
\[
\cdots z_{q(n-1)+q(n-1)-1}bx_{q(n-1)+q(n-1)}by_{q(n-1)+q(n-1)}a)\in N^{2}%
\]
and $2q(n-1)=q(n)$ proves that the statement is valid in $R$. $\square$

\bigskip

\noindent3. COHEN'S THEOREM FOR LIE NILPOTENT RINGS

\bigskip

\noindent\textbf{3.1.Proposition.} \textit{Let }$R$\textit{ be a ring with
}$\mathrm{L}_{n}$\textit{ (}$n\geq2$\textit{) and }$a,b\in R$\textit{.}

\begin{enumerate}
\item \textit{If }$P\vartriangleleft R$\textit{ is a prime ideal and }$ab\in
P$\textit{, then }$a\in P$\textit{ or }$b\in P$\textit{. In other words: all
prime ideals of }$R$\textit{\ are completely prime.}

\item \textit{The prime radical of }$R$\textit{ is the set of all nilpotent
elements:}%
\[
\mathrm{rad}(R)=\{u\in R\mid u^{k}=0\text{\textit{ for some }}k\geq1\}.
\]

\item \textit{The factor ring }$R/\mathrm{rad}(R)$\textit{ is commutative.}

\item \textit{If }$L\leq_{l}R$\textit{ is a left ideal of }$R$\textit{, then
}$L+\mathrm{rad}(R)\vartriangleleft R$\textit{ is a two sided ideal of }%
$R$\textit{. In particular, all left ideals containing }$\mathrm{rad}%
(R)$\textit{\ are two sided ideals.}
\end{enumerate}

\bigskip

\noindent\textbf{Proof.}

(1): The factor ring $S=R/P$ is a prime ring with $\mathrm{L}_{n}$ and $ab\in
P$ implies that%
\[
\overline{a}\overline{b}=(a+P)(b+P)=\overline{0}%
\]
in $S$. The application of Theorem 2.5 gives that%
\[
\underset{1.}{\underbrace{\overline{b}S\overline{b}S\overline{a}}}%
S\underset{2.}{\underbrace{\overline{b}S\overline{b}S\overline{a}}}S\ldots
S\underset{q.}{\underbrace{\overline{b}S\overline{b}S\overline{a}}}=\{0\},
\]
where $q=q(n)=2^{n-2}$. Since $S$ is prime, we deduce that $\overline
{a}=\overline{0}$ or $\overline{b}=\overline{0}$. Thus we have $a\in P$ or
$b\in P$.

(2): Let $P\vartriangleleft R$ be an arbitrary prime ideal of $R$\ and $u\in
R$ a nilpotent element with $u^{k}=0$. Now the iterated application of part
(1) of the present Proposition 3.1 gives that $u\in P$ and%
\[
\{u\in R\mid u^{k}=0\text{ for some }k\geq1\}\subseteq\mathrm{rad}(R)\text{.}%
\]
The reverse containment holds in any ring.

(3): Part (1) of Corollary 2.4 gives that $[x,y]\in R[R,R]R$ is nilpotent,
whence $[x,y]\in\mathrm{rad}(R)$ follows by part (2) of the present
Proposition 3.1.

(4) Since the sum of left ideals is a left ideal, we only have to show that
$(x+u)r\in L+\mathrm{rad}(R)$ for all $x\in L$, $u\in\mathrm{rad}(R)$ and
$r\in R$. We have%
\[
(x+u)r=rx+[x,r]+ur,
\]
where $rx\in L$, $ur\in\mathrm{rad}(R)$ and the nilpotency of $[x,r]$ gives
that $[x,r]\in\mathrm{rad}(R)$ as in part (3) above. $\square$

\bigskip

\noindent\textbf{3.2.Remark.} If $R$ is a\ simple (unitary) Lie nilpotent
ring, then the two sided ideal $\mathrm{rad}(R)$ is zero and the commutativity
of $R\cong R/\mathrm{rad}(R)$ implies that $R$ is a field.

\bigskip

We have now collected sufficient tools to obtain a Lie nilpotent version of
the following famous theorem in commutative algebra:

\bigskip

\noindent\textbf{Cohen's Theorem.} \textit{If every prime ideal of a
commutative ring }$R$\textit{ is finitely generated, then every ideal of }%
$R$\textit{ is finitely generated (i.e. }$R$\textit{ is Noetherian).}

\bigskip

In quite a number of proofs of some versions of Cohen's Theorem Zorn's Lemma
comes in handy (see, for example [S]), as is the case below:

\bigskip

\noindent\textbf{3.3.Theorem.} \textit{If the (completely) prime ideals of a
Lie nilpotent ring }$R$\textit{ are finitely generated as left ideals, then
every (left) ideal of }$R$\textit{ containing }$\mathrm{rad}(R)$\textit{ is
finitely generated as a left ideal.}

\bigskip

\noindent\textbf{Proof.} For the sake of contradiction suppose that%
\[
\mathcal{N}=\{L\subseteq R\mid L\text{ is a non-finitely generated left ideal
of }R\text{\ and }\mathrm{rad}(R)\subseteq L\}
\]
is a non-empty set. In view of part (4) of Proposition 3.1 any left ideal
$L\subseteq R$ with $\mathrm{rad}(R)\subseteq L$ is a two sided ideal of $R$.
Thus the elements of $\mathcal{N}$ are two sided ideals. A straigthforward
argument shows that the union of the (left) ideals of a chain (with respect to
the containment relation) in $\mathcal{N}$ is also an element of $\mathcal{N}%
$. The application of Zorn's lemma gives the existence of a maximal element
$P$ in $\mathcal{N}$. We claim that $P$ is a completely prime ideal of $R$.

Assume that $ab\in P$ and $a,b\in R\smallsetminus P$. The maximality of $P$
and $P\subseteq P+Rb$, $b\in P+Rb$ imply that the (left) ideal%
\[
P+Rb=R(p_{1}+r_{1}b)+\cdots+R(p_{k}+r_{k}b)
\]
is finitely generated by some elements $p_{i}+r_{i}b$, $1\leq i\leq k$ with
$p_{i}\in P$ and $r_{i}\in R$. Consider the set%
\[
K=\{x\in R\mid xb\in P\}.
\]
Clearly, $K\subseteq R$ is a left ideal of $R$\ and $a\in K$. The containment
$P\subseteq K$ follows from the fact that $P$ is a two sided ideal. The
maximality of $P$ ensures that $K$ is a finitely generated left ideal, whence
we obtain that $Kb\subseteq R$ is also a finitely generated left ideal. We
claim that%
\[
P=Rp_{1}+\cdots+Rp_{k}+Kb.
\]
Since $p_{i}\in P$ and $Kb\subseteq P$, we have $Rp_{1}+\cdots+Rp_{k}%
+Kb\subseteq P$. On the other hand, an element $p\in P\subseteq P+Rb$ can be
written as%
\[
p=s_{1}(p_{1}+r_{1}b)+\cdots+s_{k}(p_{k}+r_{k}b),
\]
whence%
\[
(s_{1}r_{1}+\cdots+s_{k}r_{k})b=p-s_{1}p_{1}-\cdots-s_{k}p_{k}\in P
\]
and $s_{1}r_{1}+\cdots+s_{k}r_{k}\in K$ follow. Thus we have%
\[
p=(s_{1}p_{1}+\cdots+s_{k}p_{k})+(s_{1}r_{1}+\cdots+s_{k}r_{k})b\in
Rp_{1}+\cdots+Rp_{k}+Kb.
\]
and $P\subseteq Rp_{1}+\cdots+Rp_{k}+Kb$.

Since $Kb$ is a finitely generated left ideal, we obtain that $P=Rp_{1}%
+\cdots+Rp_{k}+Kb$ is also a finitely generated left ideal of $R$, a
contradiction. Thus $ab\in P$ and $a,b\in R\smallsetminus P$ is impossible,
proving that $P$ is completely prime.

Now $P\in\mathcal{N}$ contradicts the condition that all completely prime
ideals are finitely generated as left ideals. It follows that $\mathcal{N}%
=\varnothing$ and our proof is complete. $\square$

\bigskip

\noindent\textbf{3.4.Remark.} Since Theorem 3.3 concerns the left ideals
containing the prime radical, it is not a full generalization of Cohen's
Theorem. On the other hand in a certain sense Theorem 3.3 is stronger than the
existing non-commutaive generalizations of Cohen's Theorem. The reason is that
prime left ideals are not used, we impose conditions only on the two sided
prime ideals.

\bigskip

\noindent\textbf{3.5.Remark.} Since in the Lie nilpotent case $\overline
{R}=R/\mathrm{rad}(R)$ is commutative (see (3) of Proposition 3.1), a direct
application of Cohen's original theorem to $\overline{R}$\ gives the following
weaker version of Theorem 3.3:

\textit{If the (completely) prime ideals and the prime radical }%
$\mathrm{rad}(R)$\textit{\ of a Lie nilpotent }$R$\textit{ are finitely
generated as left ideals, then every (left) ideal of }$R$\textit{ containing
}$\mathrm{rad}(R)$\textit{ is finitely generated as a left ideal.}

It seems that the containment of the prime radical cannot be omitted in this
direct application.

\bigskip

\noindent4. THE $n$-TH LIE CENTER

\bigskip

Let $R$ be an arbitrary ring with $1$, the $n$\textit{-th Lie center} of $R$
is defined as%
\[
\mathrm{Z}_{n}(R)=\{r\in R\mid\lbrack r,x_{1},\ldots,x_{n}]_{n+1}^{\ast
}=0\text{ for all }x_{i}\in R\text{, }1\leq i\leq n\}.
\]
The fact that $\mathrm{Z}_{n}(R)$ is a (unitary) subring of $R$ is a
consequence of $[rs,x]=[r,sx]+[s,xr]$ and%
\[
\mathrm{Z}(R)=\mathrm{Z}_{1}(R)\subseteq\mathrm{Z}_{2}(R)\subseteq
\cdots\subseteq\mathrm{Z}_{n}(R)\subseteq\mathrm{Z}_{n+1}(R)\subseteq\cdots
\]
follows from%
\[
\lbrack r,x_{1},\ldots,x_{n},x_{n+1}]_{n+2}^{\ast}=[[r,x_{1},\ldots
,x_{n}]_{n+1}^{\ast},x_{n+1}].
\]
Since%
\[
\lbrack\lbrack r,s],x_{1},\ldots,x_{n}]_{n+1}^{\ast}=[[r,s,x_{1},\ldots
x_{n-1}]_{n+1}^{\ast},x_{n}],
\]
$r\in\mathrm{Z}_{n}(R)$ implies $[r,s]\in\mathrm{Z}_{n}(R)$ for all $s\in R$,
so that $\mathrm{Z}_{n}(R)$ is a Lie ideal.

In any Lie ring $[x_{1},\ldots,x_{k},r]_{k+1}^{\ast}$ can be written as a sum
of $2^{k-1}$\ elements of the form $\pm\lbrack r,x_{\pi(1)},\ldots,x_{\pi
(k)}]_{k+1}^{\ast}$, where $\pi$ is some permutation of $\{1,2,\ldots,k\}$
(the use of the Jacobi identity and an easy induction on $k$ works). It
follows that $[x_{1},\ldots,x_{k},r,x_{k+1},\ldots,x_{n}]_{n+1}^{\ast}$ can be
written as a sum of some
\[
\pm\lbrack r,x_{\pi(1)},\ldots,x_{\pi(k)},x_{k+1},\ldots,x_{n}]_{n+1}^{\ast}.
\]
Thus $r\in\mathrm{Z}_{n}(R)$ implies that $[x_{1},\ldots,x_{k},r,x_{k+1}%
,\ldots,x_{n}]_{n+1}^{\ast}=0$ for all $x_{i}\in R$, $1\leq i\leq n$. Consider
the elements $r=E_{2,3}$ and $y_{1}=E_{3,4}$, $y_{2}=E_{1,2}$ in the
$K$-subalgebra $R=KI_{4}+R_{4}(1,1,1,1)$ of $\mathrm{M}_{4}(K)$ (see the
example in Section 2). For $x_{1},x_{2}\in R$ we have $\alpha,\beta,\gamma\in
K$ such that $[x_{1},x_{2}]=\alpha E_{1,3}+\beta E_{1,4}+\gamma E_{2,4}$. Thus%
\[
\lbrack\lbrack x_{1},x_{2}],r]=0\text{ and }[[r,y_{1}],y_{2}]=-E_{1,4}\neq0
\]
show that the implication%
\[
\lbrack\lbrack x_{1},x_{2}],r]=0\text{ for all }x_{1},x_{2}\in
R\Longrightarrow r\in\mathrm{Z}_{2}(R)
\]
(the converse of the mentioned one) is not valid.

The ring $\mathrm{Z}_{n}(R)$ obviously has the $\mathrm{L}_{n}$ property, a
much stronger statement is the following.

\bigskip

\noindent\textbf{4.1.Theorem.} \textit{Let }$n\geq1$\textit{ be an integer and
}$C\subseteq R$\textit{ a commutative submonoid of the multiplicative monoid
of }$R$\textit{. Then the subring }$S=\langle\mathrm{Z}_{n}(R)\cup C\rangle
$\textit{ of }$R$\textit{ generated by the subset }$\mathrm{Z}_{n}(R)\cup
C\subseteq R$\textit{ also has the }$\mathrm{L}_{n}$\textit{ property, i.e.}%
\[
\lbrack x_{1},x_{2},\ldots,x_{n},x_{n+1}]_{n+1}^{\ast}=0
\]
\textit{is a polynomial identity on }$S$\textit{.}

\bigskip

\noindent\textbf{Proof.} Since $cr=rc-[r,c]$ and $[r,c]\in\mathrm{Z}_{n}(R)$
for all $r\in\mathrm{Z}_{n}(R)$ and $c\in C$, we deduce that any element of
the subring $S=\langle\mathrm{Z}_{n}(R)\cup C\rangle$ can be written as%
\[
r_{1}c_{1}+\cdots+r_{t}c_{t}%
\]
with $r_{1},\ldots,r_{t}\in\mathrm{Z}_{n}(R)$ and $c_{1},\ldots,c_{t}\in C$
(notice that $1\in\mathrm{Z}_{n}(R)\cap C$).

In order to check that $[x_{1},x_{2},\ldots,x_{n},x_{n+1}]_{n+1}^{\ast}=0$ is
a polynomial identity on $S$, it is enough to consider substitutions of the
form%
\[
x_{1}=r_{1}c_{1},\ldots,x_{n}=r_{n}c_{n},x_{n+1}=r_{n+1}c_{n+1}%
\]
with $r_{i}\in\mathrm{Z}_{n}(R)$ and $c_{i}\in C$ ($1\leq i\leq n+1$).

If $1\leq k\leq n$, then we claim that%
\[
\lbrack r_{1}c_{1},\ldots,r_{n}c_{n},r_{n+1}c_{n+1}]_{n+1}^{\ast}=
\]%
\[
r_{n+1}[r_{n}[r_{n-1}[\ldots r_{k+2}[r_{k+1}[[r_{1}c_{1},\ldots,r_{k}%
c_{k}]_{k}^{\ast},c_{k+1}],\ldots],c_{n}],c_{n+1}]
\]
holds for any choice of the elements $r_{i}\in\mathrm{Z}_{n}(R)$ and $c_{i}\in
C$ ($1\leq i\leq n+1$).

Clearly, $r_{n+1}\in\mathrm{Z}_{n}(R)$ implies that%
\[
\lbrack r_{1}c_{1},\ldots,r_{n}c_{n},r_{n+1}c_{n+1}]_{n+1}^{\ast}=[[r_{1}%
c_{1},\ldots,r_{n}c_{n}]_{n}^{\ast},r_{n+1}c_{n+1}]=
\]%
\[
r_{n+1}[[r_{1}c_{1},\!\ldots\!,r_{n}c_{n}]_{n}^{\ast},\!c_{n+1}]\!+\![[r_{1}%
c_{1},\!\ldots\!,r_{n}c_{n}]_{n}^{\ast},\!r_{n+1}]c_{n+1}\!=\!r_{n+1}%
[[r_{1}c_{1},\!\ldots\!,r_{n}c_{n}]_{n}^{\ast},\!c_{n+1}]
\]
proving our claim for $k=n$.

In view of the commutativity of $C$ we have $[xc,c^{\prime}]=[x,c^{\prime}]c$
for all $x\in R$ and $c,c^{\prime}\in C$, whence%
\[
r_{n+1}[r_{n}[\!\ldots\![r_{k+1}[xc_{k},c_{k+1}],\!\ldots\!],c_{n}%
],c_{n+1}]\!=\!r_{n+1}[r_{n}[\!\ldots\![r_{k+1}[x,c_{k+1}],\!\ldots
\!],c_{n}],c_{n+1}]c_{k}%
\]
follows. Now assume that our claim holds for some $2\leq k\leq n$, then we
obtain that%
\[
\lbrack r_{1}c_{1},\ldots,r_{n}c_{n},r_{n+1}c_{n+1}]_{n+1}^{\ast}=
\]%
\[
r_{n+1}[r_{n}[\ldots\lbrack r_{k+1}[[r_{1}c_{1},\ldots,r_{k}c_{k}]_{k}^{\ast
},c_{k+1}],\ldots],c_{n}],c_{n+1}]=
\]%
\[
r_{n+1}[r_{n}[\ldots\lbrack r_{k+1}[[[r_{1}c_{1},\ldots,r_{k-1}c_{k-1}%
]_{k-1}^{\ast},r_{k}c_{k}],c_{k+1}],\ldots],c_{n}],c_{n+1}]=
\]%
\[
r_{n+1}[r_{n}[\ldots\lbrack r_{k+1}[r_{k}[[r_{1}c_{1},\ldots,r_{k-1}%
c_{k-1}]_{k-1}^{\ast},c_{k}],c_{k+1}],\ldots],c_{n}],c_{n+1}]+
\]%
\[
r_{n+1}[r_{n}[\ldots\lbrack r_{k+1}[[[r_{1}c_{1},\ldots,r_{k-1}c_{k-1}%
]_{k-1}^{\ast},r_{k}]c_{k},c_{k+1}],\ldots],c_{n}],c_{n+1}]=
\]%
\[
r_{n+1}[r_{n}[\ldots\lbrack r_{k+1}[r_{k}[[r_{1}c_{1},\ldots,r_{k-1}%
c_{k-1}]_{k-1}^{\ast},c_{k}],c_{k+1}],\ldots],c_{n}],c_{n+1}]+
\]%
\[
r_{n+1}[r_{n}[\ldots\lbrack r_{k+1}[[[r_{1}c_{1},\ldots,r_{k-1}c_{k-1}%
]_{k-1}^{\ast},r_{k}]c_{k},c_{k+1}],\ldots],c_{n}],c_{n+1}]=
\]%
\[
r_{n+1}[r_{n}[\ldots\lbrack r_{k+1}[r_{k}[[r_{1}c_{1},\ldots,r_{k-1}%
c_{k-1}]_{k-1}^{\ast},c_{k}],c_{k+1}],\ldots],c_{n}],c_{n+1}]+
\]%
\[
r_{n+1}[r_{n}[\ldots\lbrack r_{k+1}[[[r_{1}c_{1},\ldots,r_{k-1}c_{k-1}%
]_{k-1}^{\ast},r_{k}],c_{k+1}],\ldots],c_{n}],c_{n+1}]c_{k}%
\]
for all $r_{i}\in\mathrm{Z}_{n}(R)$ and $c_{i}\in C$ ($1\leq i\leq n+1$).
Since $r_{k}\in\mathrm{Z}_{n}(R)$, the substitution of $c_{k}=1\in C$ in the
above identity gives that%
\[
0=[r_{1}c_{1},\ldots,r_{k-1}c_{k-1},r_{k},r_{k+1}c_{k+1},\ldots,r_{n+1}%
c_{n+1}]_{n+1}^{\ast}=
\]%
\[
r_{n+1}[r_{n}[\ldots\lbrack r_{k+1}[r_{k}[[r_{1}c_{1},\ldots,r_{k-1}%
c_{k-1}]_{k-1}^{\ast},1],c_{k+1}],\ldots],c_{n}],c_{n+1}]+
\]%
\[
r_{n+1}[r_{n}[\ldots\lbrack r_{k+1}[[[r_{1}c_{1},\ldots,r_{k-1}c_{k-1}%
]_{k-1}^{\ast},r_{k}],c_{k+1}],\ldots],c_{n}],c_{n+1}]=
\]%
\[
r_{n+1}[r_{n}[\ldots\lbrack r_{k+1}[[[r_{1}c_{1},\ldots,r_{k-1}c_{k-1}%
]_{k-1}^{\ast},r_{k}],c_{k+1}],\ldots],c_{n}],c_{n+1}],
\]
whence%
\[
\lbrack r_{1}c_{1},\ldots,r_{n}c_{n},r_{n+1}c_{n+1}]_{n+1}^{\ast}=
\]%
\[
r_{n+1}[r_{n}[\ldots\lbrack r_{k+1}[r_{k}[[r_{1}c_{1},\ldots,r_{k-1}%
c_{k-1}]_{k-1}^{\ast},c_{k}],c_{k+1}],\ldots],c_{n}],c_{n+1}]
\]
follows. Thus the validity of our claim inherits from $k$ to $k-1$.

For $k=1$ the above claim gives that%
\[
\lbrack r_{1}c_{1},\ldots,r_{n}c_{n},r_{n+1}c_{n+1}]_{n+1}^{\ast}%
=r_{n+1}[r_{n}[\ldots\lbrack r_{2}[r_{1}c_{1},c_{2}],\ldots],c_{n}],c_{n+1}]
\]
for all $r_{i}\in\mathrm{Z}_{n}(R)$ and $c_{i}\in C$ ($1\leq i\leq n+1$). Take
$c_{1}=1\in C$ in the above identity and use $r_{1}\in\mathrm{Z}_{n}(R)$ to
derive%
\[
0=[r_{1},r_{2}c_{2},\ldots,r_{n}c_{n},r_{n+1}c_{n+1}]_{n+1}^{\ast}%
=r_{n+1}[r_{n}[\ldots\lbrack r_{2}[r_{1},c_{2}],\ldots],c_{n}],c_{n+1}],
\]
whence%
\[
\lbrack r_{1}c_{1},\ldots,r_{n}c_{n},r_{n+1}c_{n+1}]_{n+1}^{\ast}%
=r_{n+1}[r_{n}[\ldots\lbrack r_{2}[r_{1}c_{1},c_{2}],\ldots],c_{n}],c_{n+1}]=
\]%
\[
r_{n+1}[r_{n}[\ldots\lbrack r_{2}[r_{1},c_{2}],\ldots],c_{n}],c_{n+1}]c_{1}=0
\]
follows. $\square$

\bigskip

\noindent ACKNOWLEDGMENT

\noindent The authors sincerely thank the referee for many helpful comments
which helped to improve the exposition of the original submitted paper.

\bigskip

\noindent REFERENCES

\bigskip

\begin{enumerate}
\item {}[Ch] V. R. Chandran, \textit{On two analogues of Cohen's theorem,}
Pure Appl. Math.~Sci. 7(1-2) (1978), 5-10.

\item {}[Co] I. S. Cohen, \textit{Commutative rings with restricted minimum
condition,} Duke Math.~J. 17 (1950), 27-42.

\item {}[D] V. Drensky, \textit{Free Algebras and PI-Algebras,}
Springer-Verlag, New York, 2000.

\item {}[DF] V. Drensky and E. Formanek, \textit{Polynomial Identity Rings,}
Birkh\"{a}user-Verlag, Basel, 2004.

\item {}[H] F. Hansen, \textit{On one-sided prime ideals}, Pacific J.~Math.
58(1) (1975), 79-85.

\item {}[Je] S. A. Jennings, \textit{On rings whose associated Lie rings are
nilpotent,} Bull. Amer. Math. Soc. 53 (1947), 593-597.

\item {}[Jo] P. Jothilingam, \textit{Cohen's theorem and Eakin-Nagata theorem
revisited}, Comm..~Algebra 28(10) (2000), 4861-4866.

\item {}[Ka] I. Kaplansky, \textit{Elementary divisors and modules},
Trans.~Amer.~Math. Soc. 66 (1949), 464-491.

\item {}[Ko] K. Koh, \textit{On prime one-sided ideals}, Canad.~Math.~Bull.
14(2) (1971), 259-260.

\item {}[Kr] G. Krause, On fully left bounded left noetherian rings, J.
Algebra 23 (1972), 88--99.

\item {}[L] C.-P. Lu, \textit{Spectra of modules,} Comm.~Algebra 23(10)
(1995), 3741-3752.

\item {}[Na] A. R. Naghipour, \textit{A simple proof of Cohen's theorem},
Amer.~Math. Monthly 112(9) (2005), 825-826.

\item {}[Ni] I. Nishitani, \textit{A Cohen-type theorem for Artinian modules},
Arch.~Math. 87(3) (2006), 206-210.

\item {}[Re] M. L. Reyes, \textit{Noncommutative generalizations of theorems
of Cohen and Kaplansky,} Algebr.~Represent.~Theory 15(5) (2012), 933-975.

\item {} [Ro] L. H. Rowen, \textit{Polynomial Identities in Ring Theory,}
Academic Press, New York, 1980.

\item {}[S] R. Y. Sharp, \textit{Steps in Commutative Algebra}, Second
edition, Cambridge University Press, Cambridge, 2000.

\item {}[Sz] J. Szigeti, \textit{New determinants and the Cayley-Hamilton
theorem for matrices over Lie nilpotent rings}, Proc. Amer. Math. Soc. 125(8)
(1997), 2245-2254.\ 

\item {}[SzvW] J. Szigeti and L. van Wyk, \textit{Determinants for }$n\times
n$\textit{ matrices and the symmetric Newton formula in the }$3\times
3$\textit{ case }, Linear and Multilinear Algebra, Published online: 16 Jul
2013, DOI:10.1080/03081087.2013.806919
\end{enumerate}

\end{document}